\documentclass[a4paper,12pt,leqno]{amsart}

\usepackage[T1]{fontenc}
\usepackage{mathrsfs}
\usepackage{latexsym}
\usepackage{epsfig}
\usepackage{mathtools}
\usepackage{amsmath,amsfonts,amsthm,amssymb,amscd}
\usepackage{color}
\usepackage{verbatim}
\usepackage{hyperref}
\usepackage{dsfont}
\usepackage{bbm}
\usepackage{tcolorbox}
\usepackage{pstricks}
\usepackage{float}
\usepackage{tikz}
\usepackage{pgfplots}
\usepackage{tikz-3dplot}
\pgfplotsset{compat=newest}
\usepackage{subcaption}
\usetikzlibrary{patterns.meta,perspective}
\pgfplotsset{compat=1.4}
\usetikzlibrary{decorations.pathreplacing}
\usepackage{tkz-euclide}
\newcommand\be{\begin{equation}}
\newcommand\ee{\end{equation}}
\newcommand\benn{\begin{equation*}}
\newcommand\eenn{\end{equation*}}
\newcommand\bc{\begin{center}}
\newcommand\ec{\end{center}}
\newcommand\bea{\begin{eqnarray}}
\newcommand\eea{\end{eqnarray}}
\newcommand\beann{\begin{eqnarray*}}
\newcommand\eeann{\end{eqnarray*}}
\newcommand\bi{\begin{itemize}}
\newcommand\ei{\end{itemize}}
\newcommand\ben{\begin{enumerate}}
\newcommand\een{\end{enumerate}}
\newcommand\bml{\begin{multline}}
\newcommand\eml{\end{multline}}
\newcommand\bmln{\begin{multline*}}
\newcommand\emln{\end{multline*}}

\newtheorem{thm}{Theorem}[section]

\newtheorem{cor}[thm]{Corollary}
\newtheorem{lem}[thm]{Lemma}

\newtheorem{rek}[thm]{Remark}

\newtheorem{cla}[thm]{Claim}

\newcommand{\avec}[1]{\langle #1 \rangle}
\newcommand{\h}[1]{\widehat{#1}}
\newcommand{\lla}[1]{\left \{  #1  \right \}}
\newcommand{\pa}[1]{\left (  #1 \right )}
\newcommand{\f}[1]{\vec{#1}}
\newcommand{\abs}[1]{\left |  #1  \right |}
\newcommand{\nor}[1]{\left \|  #1  \right \|}
\newcommand{\corch}[1]{\left [  #1  \right ] }

\newcommand{\R}{\ensuremath{\mathbb{R}}}

\newcommand{\N}{\mathbb{N}}

\newcommand{\dps}{\displaystyle}

\newcommand{\df}{\mathrm{d}}
\newcommand{\ve}{\varepsilon}
\newcommand{\Leb}{\mathcal{L}}
\newcommand{\Ha}{\mathcal{H}}
 
\newcommand{\tri}{\bigtriangleup}

\newcommand{\floor}[1]{\left\lfloor #1 \right\rfloor}
\newcommand{\qedd}{\hfill$\square$}

\begin{document}

\definecolor{qqqqff}{rgb}{0.,0.,1.}
\definecolor{ffqqqq}{rgb}{1.,0.,0.}

\title{A Mattila-Sj\"olin theorem for simplices in low dimensions}

 \author{Eyvindur Ari Palsson}
 \email{\textcolor{blue}{\href{mailto:palsson@vt.edu} {palsson@vt.edu}}}
 \address{Department of Mathematics, Virginia Tech, Blacksburg, VA 24061}

 \author{Francisco Romero Acosta}
 \email{\textcolor{blue}{\href{mailto:jfromero@vt.edu} {jfromero@vt.edu}}}
 \address{Department of Mathematics, Virginia Tech, Blacksburg, VA 24061}
 
\begin{abstract}
In this paper we show that if a compact set $E \subset \mathbb{R}^d$, $d \geq 3$, has Hausdorff dimension greater than $\frac{(4k-1)}{4k}d+\frac{1}{4}$ when $3 \leq d<\frac{k(k+3)}{(k-1)}$ or $d- \frac{1}{k-1}$ when $\frac{k(k+3)}{(k-1)} \leq d$, then the set of congruence class of simplices with vertices in $E$ has nonempty interior. By set of congruence class of simplices with vertices in $E$ we mean
\begin{center}
  $\Delta_{k}(E) = \left \{ \f{t} = \pa{t_{ij}} : |x_i-x_j|=t_{ij} ; \   x_i,x_j \in E ; \  0\leq i < j \leq k \right \} \subset \mathbb{R}^{\frac{k(k+1)}{2}}$
\end{center}
  where $2 \leq k <d$. This result improves the previous best results in the sense that we now can obtain a Hausdorff dimension threshold which allow us to guarantee that the set of congruence class of triangles formed by triples of points of $E$ has nonempty interior when $d=3$ as well as extending to all simplices. The present work can be thought of as an extension of the Mattila-Sj\"olin theorem which establishes a non-empty interior for the distance set instead of the set of congruence classes of simplices. 
\end{abstract}
\maketitle

\section{Introduction}
Falconer's distance conjecture states that if the Hausdorff dimension of a compact set $E \subset \R^d$, $d \geq 2$, is greater than $\frac{d}{2}$, then its distance set $\tri(E)= \lla{|x-y|; \ x,y \in E}$ has positive Lebesgue measure. Falconer \cite{Falconer} not only stated this conjecture, but he also proved that if $\dim_{\Ha}(E)>\frac{d+1}{2}$, then $\tri(E)$ has positive Lebesgue measure. Falconer's conjecture remains open, but much work has been done towards it \cite{Bourgain,Erdogan,Liu,Wolff1}. For instance, when $d=2$ Bourgain \cite{Bourgain} showed that if $\dim_{\Ha}(E)>\frac{13}{9}$, then $\Leb^1 \pa{\tri(E)}>0$, later Wolff \cite{Wolff1} improved the threshold to $\frac{4}{3}$. In dimension $3$ and higher, Erdogan \cite{Erdogan} obtained the threshold $\frac{d}{2} + \frac{1}{3}$. Currently, the best known result when $d=2$ is due to Guth, Iosevich, Ou and Wang \cite{GuthIosevichOuWang}, who showed that if $\dim_{\Ha}(E)>\frac{5}{4}$, then $\tri(E)$ has positive Lebesgue measure. When $d=3$, Du, Guth, Ou, Wang, Wilson and Zhang \cite{DuGuthOuWangWilsonZhang} showed that if $\dim_{\Ha}(E)>\frac{9}{5}$, then $\Leb^1\pa{\tri(E)}>0$ while for higher dimensions, that is for $d\geq 4$, Du and Zhang \cite{DuZhang} improved the threshold to $\frac{d^2}{2d-1}$. Further, when restricting $d \geq 4$ to even integers, Du, Iosevich, Ou, Wang and Zhang \cite{DuIosevichOuWangZhang} showed that $\dim_{\Ha}(E)>\frac{d}{2} + \frac{1}{4}$ is enough to guarantee $\Leb^1 \pa{\tri(E)} >0$. Most recently, Du, Ou, Ren and Zhang \cite{DuOuRenZhang} broke the $\frac{d}{2} + \frac{1}{4}$ barrier when $d\geq 3$ and obtained the threshold $\frac{d}{2}+\frac{1}{4}-\frac{1}{8d+4}$. Dimension versions of Falconer's distance conjecture have also been studied, for details see the work done by Shmerkin and Wang \cite{ShmerkinWang} and the references therein.

A classic result due to H. Steinhaus \cite{Steinhaus} states that if a set $E \subset \R^d$, $d \geq 1$, has positive Lebesgue measure, then the set $E-E = \lla{ x-y : x,y \in E} \subset \R^d$ contains a neighborhood of the origin. Likewise, in the context of distance set problem we might wonder the following: For a given compact set $E \subset \R^d$, how large does its Hausdorff dimension need to be to guarantee that its distance set contains an interval. Mattila and Sj\"olin \cite{MattilaSjolin} proved that if $\dim_{\Ha}(E)>\frac{d+1}{2}$, then its distance set $\tri(E)$ has nonempty interior. The proof of this result can be found in \cite{Mattila1}. Iosevich, Mourgoglou and Taylor \cite{IosevichMourgoglouTaylor} showed that if $\dim_{\Ha}(E)>\frac{d+1}{2}$, then the set $\tri (E) = \lla{ \nor{x-y}_{B}: x,y \in E }$ contains an interval. Here, $ \nor{ \cdot }_{B}$ is a metric induced by the norm defined by a bounded convex body $B$ with non-vanishing curvature. Koh, Pham and Shen \cite{KohPhamShen} obtain slight improvements to these thresholds in the case of product sets, that is, sets of the form $E = A \times A \times \hdots \times A = A^d$, where $A \subset \R$ is compact. Greenleaf, Iosevich and Taylor \cite{GreenleafIosevichTaylor1} extend the Mattila-Sj\"olin theorem to more general 2-point configuration sets, that is, they obtained a Hausdorff dimensional threshold for which they can guarantee that the set $\Delta_{\Phi}(E) = \lla{ \Phi(x,y): x,y \in E }$ of $\Phi-$configurations in $E$ has nonempty interior in $\R^m$, $m \geq 1$ where $\Phi: \R^d \times \R^d \to \R^m$ is a smooth function satisfying some additional conditions. Note that by considering $\Phi(x,y)=|x-y|$, we have that $\Delta_{\Phi}(E) = \Delta(E)$, in this case the authors in \cite{GreenleafIosevichTaylor1} (Corollary $1.2$) obtained a Hausdorff dimensional threshold that coincides with the threshold obtained in \cite{MattilaSjolin}.

More generally, one can study analogues of the Falconer distance problem and the Mattila-Sj\"olin theorem for $(k+1)$-point configurations. The particular $(k+1)$-point configurations we study in this paper are simplices. We define the set of congruence classes of simplices with vertices in $E$ as the set
\bc
  $\dps \Delta_{k}(E) = \left \{ \f{t}=(t_{ij}): |x_i-x_j|=t_{ij}; \   x_i, x_j, \in E, \ 0\leq i <j \leq k \right \}$.
\ec
where $\f{t}=(t_{ij})$ is an element of $\R^{\frac{k(k+1)}{2}}$. We note that with this notation the $(k+1)$-points that form the configuration are labeled by $x_0, x_1, \ldots, x_k$ and we call the corresponding simplex the $k$-simplex. Greenleaf and Iosevich \cite{GreenleafIosevich} introduced these types of questions, when they studied triangles in the plane, which corresponds to the case of $k=2$ and $d=2$. In their paper they showed that if $\dim_{\Ha}(E)>\frac{7}{4}$ then $\Leb^3\pa{\Delta_{2}(E)}>0$, which is an analogue of the Falconer distance problem. This was extended by Erdogan, Hart and Iosevich \cite{ErdoganHartIosevich}, who obtained the threshold $\frac{d+k+1}{2}$ for the $k$-simplex in $\mathbb{R}^d$, and Grafakos, Greenleaf, Iosevich and the first named author in \cite{GrafakosGreenleafIosevichPalsson}, where they obtained the threshold $d - \frac{d-1}{2k}$, which is an improvement in lower dimensions. Further, \cite{GrafakosGreenleafIosevichPalsson} included a fairly general mechanism that worked for a host of other $(k+1)$-point configurations. The current best results on these problems in low dimensions are due to Greenleaf, Iosevich, Liu and the first named author \cite{GreenleafIosevichLiuPalsson2} where they obtain the threshold $\frac{2d^2}{3d-1}$ while in high enough dimensions the best threshold is $\frac{d+k}{2}$ due to Iosevich, Pham, Pham and Shen \cite{IosevichPhamShen}.

The first result on non-empty interior for $(k+1)$-point configurations, in the spirit of the Mattila-Sj\"olin theorem for distances, is due to Bennett, Iosevich and Taylor \cite{BennettIosevichTaylor} for configurations called chains. This was further extended to configurations called trees by Iosevich and Taylor \cite{IosevichTaylor}. More recently Greenleaf, Iosevich and Taylor \cite{GreenleafIosevichTaylor2} extended their techniques from \cite{GreenleafIosevichTaylor1} to be able to handle somewhat general classes of $(k+1)$-point configurations. Their techniques did not yield results for simplices which motivated the question of how large the Hausdorff dimension of a compact set $E \subset \R^d$ needs to be to guarantee that the set $\Delta_{k}(E)$ has nonempty interior, which is the question we study in this paper.

 In our previous work \cite{PalssonRomero}, given a compact subset $E \subset \R^d$ we used the classical rule \emph{side-angle-side} to define a set of congruence classes of triangles formed by triples of points of $E$, that is, 
\benn 
\Delta_{\text{tri}}(E) = \left \{ (t,r, \alpha) : |x-z|=t, |y-z|=r \text{\ and \ } \alpha= \alpha(x,z,y), \  x,y,z \in E \right \},
\eenn
where $\alpha (x,z,y)$ denotes the angle formed by $x$, $y$ and $z$, centered at $z$. Moreover, we showed that if a compact set $E \subset \R^d$, $d \geq 4$, has Hausdorff dimension greater than $\frac{2d}{3}+1$, then the set $\Delta_{\text{tri}}(E)$ has nonempty interior. However, this result yields nothing when $d=3$. In \cite{GreenleafIosevichTaylor2}, Greenleaf, Iosevich and Taylor applied their main result to a wide variety of $(k+1)-$point configurations. Some applications to $3-$point configurations included in \cite{GreenleafIosevichTaylor2} were area of triangles in $\R^2$, volumes of pinned parallelepipeds in $\R^3$ and ratios of pinned distances in $\R^2$ and $\R^3$. However, their results were not directly applicable to the triangle problem, regardless of how the triangles were encoded. In what was a very recent preprint when the first version of this paper was written, and what is now a published paper, Greenleaf, Iosevich and Taylor \cite{GreenleafIosevichTaylor3} refined their approach and gave an alternative proof of \cite{PalssonRomero}. Despite the sophisticated method developed by the authors in \cite{GreenleafIosevichTaylor3} it still provides a trivial threshold for the triangle problem when $d=3$.  

Similarly, we can also use the rule \emph{side-side-side} to define a set of congruence classes of triangles formed by triples of points of $E$, namely $\Delta_{2}(E)$. It is clear that there is a bijection between $\Delta_{\text{tri}}(E)$ and $\Delta_{2}(E)$. Thus, a Hausdorff dimensional threshold that guarantees that $\Delta_2(E)$ has nonempty interior, will also guarantee that $\Delta_{\text{tri}}(E)$ has nonempty interior. In the present work, we improve the result given in \cite{PalssonRomero} in the sense that we now can obtain a nontrivial threshold when $d=3$ for which we can guarantee that $\Delta_2(E)$ has nonempty interior. Moreover, we provide an analogous result for simplices in higher dimensions. We now state our main result.
\begin{thm} \label{thm:maunthm1}
Let $E \subset \mathbb{R}^d$ be a compact set, $d \geq 3$. The set of equivalence classes of $k$-simplices, 
\bc
  $\dps \Delta_{k}(E) = \left \{ \f{t}=(t_{ij}): |x_i-x_j|=t_{ij}; \   x_i, x_j, \in E, \ 0\leq i <j \leq k \right \}$,
\ec
$2 \leq k<d$, has nonempty interior if 
\benn 
    \dim_{\Ha}(E) > \left\{\begin{matrix} \frac{(4k-1)}{4k}d+\frac{1}{4} & , & 3 \leq d<\frac{k(k+3)}{(k-1)} \\ d- \frac{1}{k-1} & , &  \frac{k(k+3)}{(k-1)} \leq d \end{matrix}\right. .
\eenn 
\end{thm}

A straightforward consequence of Theorem \ref{thm:maunthm1} is the following corollary:
\begin{cor} \label{thm:maunthm2}
Let $E \subset \mathbb{R}^d$ be a compact set, $d \geq 3$. The set of congruence classes of triangles formed by triples of points of $E$,
\benn
 \Delta_{2}(E) = \left \{ (t_{01},t_{02},t_{12}) : |x_0-x_1|=t_{01}, |x_0-x_2|=t_{02}, |x_1-x_2|=t_{12} ; \  x_0, x_1, x_2 \in E \right \},
\eenn
has nonempty interior if
\benn 
    \dim_{\Ha}(E) > \left\{\begin{matrix} \frac{7d}{8}+\frac{1}{4} & , & 3 \leq d< 10 \\ d- 1 & , & 10 \leq d \end{matrix}\right. .
\eenn 

\end{cor}

As we indicated above, the main result given by Greenleaf, Iosevich and Taylor \cite{GreenleafIosevichTaylor2} allows one to obtain Hausdorff dimensional thresholds which guarantee that many general $(k+1)-$point configuration sets have nonempty interior. As they can recover the results from \cite{PalssonRomero} on triangles one can speculate whether their techniques will work in general for $k$-simplices, $k\geq 3$. The authors comment on this in \cite{GreenleafIosevichTaylor2} and note that the conditions they need for their tools appear to fail for $k\geq 3$. Thus, as the authors point out, it would be interesting to see if their microlocal analysis tools could be developed further to handle higher dimensional simplices. We thus emphasize that the novelty of this paper is not only to obtain non-trivial thresholds for triangles when $d=3$, but also obtaining non-trivial thresholds for all $k$-simplices.

As a final comment we note that these analogues of the Falconer distance problem and the Mattila-Sj\"olin theorem for simplices tell us something about the abundance of simplices in large enough sets. One can ask even stronger questions, such as if you pick a favorite configuration, e.g. an equilateral triangle, then is that configuration guaranteed to exist in a large enough set? Iosevich and Liu \cite{IosevichLiu} answered this positively in $\mathbb{R}^d$ with $d\geq 4$ and their work played an important role in our previous paper \cite{PalssonRomero}. Recently Iosevich and Magyar \cite{IosevichMagyar} lowered the dimensional threshold in \cite{IosevichLiu} to $d\geq 3$ and extended the work to simplices. Techniques from that paper play an important role in the current paper.

\subsection{Sharpness:} The threshold $\frac{d}{2}$ in Falconer's distance conjecture comes from an explicit construction by Falconer \cite{Falconer} of sets of lower dimension with distance sets of zero Lebesgue measure. For triangles a similar sharpness example exists in the plane that establishes the threshold $\frac{3}{2}$. The example is due to Erdogan and Iosevich but appeared first in print in \cite{GreenleafIosevichLiuPalsson2}. For triangles in higher dimensions and higher dimensional simplices no similar sharpness examples exist. However, since the existence of many triangles or simplices implies the existences of many distances, we have a trivial sharpness example of $\frac{d}{2}$ for these problems too. For Mattila-Sj\"olin type theorems the same sharpness examples apply and no stronger sharpness examples exist, neither for distances nor simplices in general. Thus, in the triangle case, for dimension $3$ and higher there is a gap between the dimensional threshold $\frac{7d}{8} + \frac{1}{4}$ obtained in Corollary \ref{thm:maunthm2} and the trivial threshold, $\frac{d}{2}$. In \cite{PalssonRomero}, we obtained a better dimensional threshold for dimension $4$ and higher, however a gap with respect to the trivial exponent still persists. For higher dimensional simplices we only have the results in Theorem \ref{thm:maunthm1}. There the gap between the trivial threshold and what we prove is closest in low dimensions and gets progressively worse the higher the dimension. It would be interesting to see if we could establish results analogous to in our first paper \cite{PalssonRomero} for higher dimensional simplices. \\

In the finite field setting, unlike the Euclidean one, there are results for which the thresholds between analogues of Falconer type results and analogues of Mattila-Sj\"olin type results may differ. Murphy, Petridis, Pham, Rudnev and Stevens \cite{MurphyPetridisPhamRudnevStevens} proved that if $E \subset \mathbb{F}^2_q$, $q$ prime, and $|E| \geq Cq^{\frac{5}{4}}$, for some $C>0$, then its distance set $\Delta(E)$ contains a positive proportion of $\mathbb{F}_q$. This is a result on what is known as the Erd\H{o}s-Falconer distance problem in finite fields and is an analogue of the Falconer distance problem in the Euclidean setting. Murphy and Petridis \cite{MurphyPetridis} proved that if $E \subset \mathbb{F}^2_q$ and $|E| \approx q^{\frac{4}{3}}$, then one cannot, in general, expect the distance set to contain all of $\mathbb{F}_q$. The problem of showing that the distance set contains all of $\mathbb{F}_q$ is often viewed as an analogue of the Mattila-Sj\"olin theorem in the Euclidean setting. Despite that, Iosevich and Rudnev \cite{IosevichRudnev} proved that if $E \subset \mathbb{F}^d_q$, $d \geq 2$, such that $|E| \gtrsim q^{\frac{d+1}{2}}$, then its distance set $ \Delta(E)$ contains all of $\mathbb{F}_q$ which establishes a Mattila-Sj\"{o}lin type result at a higher threshold. Understanding, whether the dimensional thresholds for the Falconer type theorems and the Mattila-Sj\"{o}lin type results in the Euclidean setting should be the same or different, is a major open problem.

\subsection{Overview of result:} Let $\mu$ be a Frostman probability measure supported on a compact subset $E \subset \R^d$, $d \geq 3$. This measure essentially encodes the dimension of the set we start with, see e.g. \cite{Mattila1} for the theory of Frostman measures. Given $\ve>0$, let $\psi_{\ve}(x)=\ve^{-d} \psi \pa{\frac{x}{\ve}} \geq 0$, where $\psi \geq 0$ is a smooth function such that its Fourier transform, $\h{\psi}$, is a  smooth compactly supported cut-off function, satisfying $\h{\psi}(0)=1$ and $0 \leq \h{\psi} \leq 1$. Let $\mu_{\ve} := \mu \ast \psi_{\ve}$, and note that $\nor{\mu_{\ve}}_{\infty} \lesssim \ve^{s-d} $. Furthermore, note that $\mu_{\ve} \to \mu$ weakly, as $\ve \to 0$, see e.g. \cite{Mattila1} for details. The proof of our main theorem will proceed as follows:
\begin{itemize}
    \item {\bf Step 1:} In order to ensure the non-degeneracy of simplices, we will show that it is possible to extract enough suitable subsets of $E$, each of positive $\mu-$measure, such that we can form non-degenerate simplices. This is established in Lemma \ref{PHprinc2} with a construction illustrated in Figure \ref{fig:cubes}.
    
    \item {\bf Step 2:} We will define the measure $\delta(\mu)$ on the set $\Delta_k(E)$ as the image of $\mu \times \hdots \times \mu$ ($k+1$ times) under the map $(x_0,\hdots,x_k) \mapsto \f{v}_k (x_0, \hdots, x_k)$, where $\f{v}_k (x_0, \hdots, x_k) \in \R^{\frac{k(k+1)}{2}}$ denotes the vector with entries $|x_i-x_j|$, $0 \leq i < j \leq k$, listed in a lexicographic order. Since $\mu$ is a probability measure we automatically have that $\delta(\mu)$ is also a probability measure. We can do the same procedure for the smoothed out versions $\mu_{\ve}$ and obtain a measure $\delta(\mu_{\ve})$, which from the definition of $\delta(.)$ will converge weakly to $\delta(\mu)$ as long as the $\mu_{\ve}$ converge weakly to $\mu$. In our previous work \cite{PalssonRomero} we also created such push-forward measures, building on the original approach from \cite{IosevichMourgoglouTaylor}. Unlike these previous papers, where the next step was to estimate the size of $\delta(\mu_{\ve})$ through a main term and an error term, we now proceed differently.
    
    \item {\bf Step 3:} We show that the density of the measure $\delta(\mu)$ is continuous through a Cauchy sequence argument. This is the first time that we are aware of that such an argument has been used to establish a Mattila-Sj\"{o}lin type theorem. To achieve this, we build on the techniques developed by Iosevich and Magyar \cite{IosevichMagyar} to obtain Lemma \ref{lemma1}. Adapting their techniques to our setting is one of the main technical contributions of this paper. Although technical in nature, Lemma \ref{lemma1} simply shows the Cauchy sequence nature of $\delta(\mu_{\ve})$ as one varies $\ve$, as long as we are above our dimensional thresholds. By using Lemma \ref{lemma1}, we will show that $\delta(\mu_{\ve})$ converges uniformly as $\ve \to 0$. Since $\delta(\mu_{\ve})$ is continuous for every $\ve>0$, the limit as $\ve \to 0$ must be a continuous function. Thus, by the uniqueness of the limit we have that the density of $\delta(\mu)$ is continuous.
\end{itemize}

\subsection{Acknowledgement: } The first listed author was supported in part by Simons Foundation Grant no. 360560. We thank Alex Iosevich for suggesting looking at \cite{IosevichMagyar} and fruitful discussions about the problem. Finally, we thank an anonymous referee for the many suggestions that significantly improved the exposition of the paper.
 
\section{Proof of Theorem \ref{thm:maunthm1}}

\subsection{Step 1:} In \cite{PalssonRomero}, given a compact subset $E \subset \R^d$ we had to ensure that the triangles formed by triples of points of $E$ were non-degenerate. This was accomplished by showing that if $E$ is large enough, then it is possible to extract three suitable subsets of $E$ apart from each other. More precisely, we used the following lemma: 
\begin{lem}\label{PHprinc} (\cite{PalssonRomero}, Lemma $2.1$).
  Let $\mu$ be a Frostman probability measure on $E \subset \R^d$, $d \geq 3$, with Hausdorff dimension greater than $\frac{2}{3} d+1$, then there are positive constants $c_1$, $c_2$, and $E_1$, $E_2$, $E_3$ subsets of $E$, such that
  \begin{itemize}
    \item [(i)] $\dps \mu(E_i) \geq c_1 >0$, for $i=1,2,3$.
    \item [(ii)] $\dps \underset{1 \leq k \leq d}{\max} \lla{ \inf \lla{ | x_k - y_k|: x \in E_i, y \in E_j \text{ \ and \ } i \neq j } } \geq c_2 >0 $, for $i,j=1,2,3 $.
  \end{itemize}
\end{lem}

\begin{rek}
The threshold given in Lemma \ref{PHprinc} was stated conveniently to match the threshold given in the main result of \cite{PalssonRomero}. In fact, the best threshold that can be obtained from the proof of Lemma \ref{PHprinc} is $\frac{d+1}{2}$. 
\end{rek}

Likewise, to ensure that the simplices under consideration are non-degenerate, we must show that it is possible to extract $k$ suitable subsets of $E$ disjoint from each other. Due to $2<k \leq d-1$, it is more than sufficient to show that we can extract $d+1$ suitable subsets apart from each other. Thus we have the following:
\begin{lem}\label{PHprinc2}
  Let $\mu$ be a Frostman probability measure on $E \subset \R^d$, with Hausdorff dimension greater than $ \frac{(d-1)\log_2(3)+ \log_2(d)+2}{1+\log_2(3)}$, then there are positive constants $\lla{c_i}^{d+1}_{i=1}$ and a collection $\lla{E_i}^{d+1}_{i=1}$ of subsets of $E$, such that
  \begin{itemize}
    \item [(i)] $\dps \mu(E_i) \geq c_i >0$, for $i \in \lla{1,2\hdots, d+1}$.
    \item [(ii)] $\dps \underset{1 \leq k \leq d}{\max} \lla{ \inf \lla{ | x_k - y_k|: x \in E_i, y \in E_j \text{ \ and \ } i \neq j } } \geq c >0 $, for all $i,j$ and $c>0$.
  \end{itemize}
\end{lem}
The proof of Lemma \ref{PHprinc2} can be found in Section $5$.

Note that when $k=2$ we have that the threshold obtained by Corollary \ref{thm:maunthm2} is $ \frac{7}{8}d+\frac{1}{4}$, which is greater than $ \frac{(d-1)\log_2(3)+ \log_2(d)+2}{1+\log_2(3)}$. Due to $\frac{(4k-1)d}{4k} + \frac{1}{4} $ is increasing with respect to $k$, then $\frac{(4k-1)d}{4k} + \frac{1}{4} > \frac{(d-1)\log_2(3)+ \log_2(d)+2}{1+\log_2(3)}$ for all $2 \leq k < d < \frac{k(k+3)}{k-1}$, and $d- \frac{1}{k-1}> \frac{(d-1)\log_2(3)+ \log_2(d)+2}{1+\log_2(3)}$ for all $d \geq \frac{k(k+3)}{k-1}$. Therefore, we can guarantee the non-degeneracy of the simplices.

\begin{figure}[h!]
    \centering
    \includegraphics[scale=0.33]{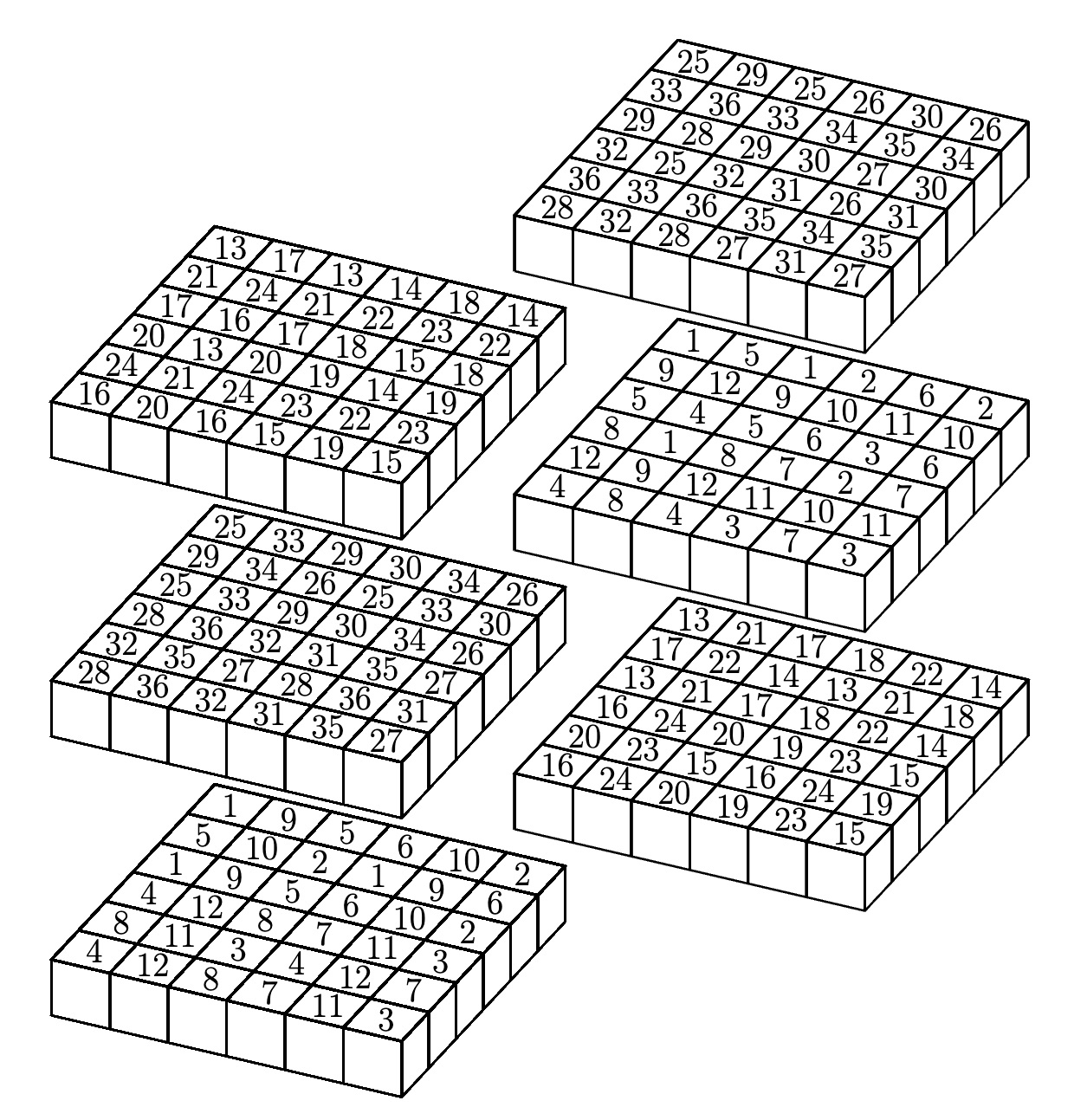}
	\caption{Due to the number of collections is greater than the number of cubes in each collection, it is simple to subdivide the unit cube into collections in a way that no more than three cubes in the same collection are 'co-planar'. For instance, when $d=3$ we have $36$ collections of $6$ cubes each. By Lemma \ref{PHprinc2} we can guarantee the existence of a collection with $4$ cubes with positive measure. The picture above shows a way in which we can subdivided the unit cube.}
	\label{fig:cubes}
\end{figure}

\begin{rek} Note the following:
\begin{itemize}
    \item The Hausdorff dimensional threshold given in Lemma \ref{PHprinc2} does not depend on $k$. This is due to the overestimation on the number of subsets that satisfy the conditions given in the Lemma. It is possible to increase the number of subsets, but this will require the set $E$ to be larger.
    \item The proof of Lemma \ref{lemma1} is given in terms of $E$ and $\mu$, but it is clear that the proof still holds for $E_i$ and $\mu_i$, $i=0,\hdots,d+1$. Where $E_i$ are the subsets that can be obtained by using Lemma \ref{PHprinc2} and $\mu_i$ are the restrictions of $\mu$ to the sets $E_i$ respectively.
    \item Note that Lemma \ref{PHprinc2} (respectively Lemma \ref{PHprinc}) ensures that the length of the edges of simplices (respectively side lengths of triangles) under consideration in Theorem \ref{thm:maunthm1} (respectively Corollary \ref{thm:maunthm2}) can be bounded above and below by positive constants. Moreover, the non-degeneracy of simplices is also guaranteed by the construction given in the proof of Lemma \ref{PHprinc2} (see Figure \ref{fig:cubes}).
\end{itemize}
\end{rek}

\subsection{Step 2:} Consider a Frostman probability measure $\mu$ supported on $E$. For any continuous function $\varphi$ on $\R^{\frac{k(k+1)}{2}}$, we defined a measure on $\dps \Delta_{k}(E)$ as follows
\benn
     \int \varphi \pa{\f{t}} \df \delta (\mu)\pa{\f{t}} = \int \hdots \int \varphi(\f{v}_k (x_0, \hdots, x_k)) \df \mu(x_0) \hdots \df \mu(x_k).
\eenn
In other words, $\delta(\mu)$ is the image of $\mu \times \hdots \times \mu$ ($k+1$ times) under the map $(x_0,\hdots,x_k) \to \f{v}_k (x_0, \hdots, v_k)$. Where $\f{v}_k (x_0, \hdots, x_k) \in \R^{\frac{k(k+1)}{2}}$ denotes the vector with entries $|x_i-x_j|$, $0 \leq i < j \leq k$, listed in a lexicographic order. Furthermore, note that for a smooth compactly supported function $f$, we have that $\delta(f)$ is also a function given by
\begin{multline}\label{deltafs}
 \delta(f)\pa{\f{t}} = \pa{\prod^{k}_{p=1} C_{F_p}} \int \hdots \int f(x) f(x+x_1) \dots  f(x+x_k) \\ \df \sigma^{d-k}_{r_k}\pa{ x_k-\sum^{k-1}_{n=1} m_{nk}x_n} \hdots \df \sigma^{d-2}_{r_2}(x_2-m_{12}x_1)  \df \sigma^{d-1}_{r_1}(x_1) \df x.
\end{multline}
Where $\df \sigma^{d-i}_{r_i}$ denotes the surface area measures over the $d-i$-dimensional sphere of radius $r_i$, and $m_{ij}$ are some positive real numbers. Furthermore, $C_{F_p} = 2^{-p} \det (A_p)^{-1/2}$, $A_p = \pa{a_{ij}}$ is a $p \times p$ matrix, where $a_{ij} = \avec{(x_p-x_{i-1}),(x_p-x_{j-1})} $, for $ 1 \leq i,j \leq p $. note that  
\be \label{eqn:mtdes}
    \dps a_{ij} = \left\{\begin{matrix} 0 & \text{, if} & (x_p-x_{i-1}) \perp (x_p-x_{j-1}) \\ t^2_{(i-1)p} & \text{, if} & i=j \\ \frac{1}{2} \pa{t^2_{(i-1)p}+t^2_{(j-1)p}-t^2_{(i-1)(j-1)}} & , \ \text{else} &  \end{matrix}\right. 
\ee
To show (\ref{deltafs}) note that  
\begin{multline*}
    \int g\pa{\f{t}} \df \delta (f)\pa{\f{t}} = \int \hdots \int g(\f{v}_k (x_0, \hdots, x_k)) f(x_0) f(x_1) \hdots f(x_k) \df x_0 \df x_1 \hdots \df x_k,
\end{multline*}
for any continuous function $g$ with compact support. By a simple change of variables we obtain 
\begin{multline*}
\int g\pa{\f{t}} \df \delta (f)\pa{\f{t}} = \int \hdots \int g(|x_1|,\hdots,|x_k|,|x_1-x_2|,|x_1-x_3|, \hdots,|x_{k-1}-x_k|) \\
f(x) f(x+x_1) \hdots f(x+x_k) \df x \df x_1 \hdots \df x_k.
\end{multline*}
Due to our previous change of variables, lets from now until the end of the proof of our main result denote $x_0=0$. Consider the polynomials $P_{ij}(w)=|w-x_i|^2-t^2_{ij}$, and let $F_j = \lla{P_{ij},i<j}$, for $0 \leq i <j \leq k$, then the right hand side of the latter equation is equal to
\begin{multline*}
    \int \hdots \int g\pa{\f{t}} \pa{ \int \hdots \int f(x) f(x+x_1) \hdots f(x+x_2) \df \omega_{F_1}(x_1) \hdots \df \omega_{F_k}(x_k)  \df x } \df \f{t},
\end{multline*}
where $\df \omega_{F_j}$ is the measure supported on the algebraic set $S_{F_j} = \left \{ x \in \R^d: P_{ij}(x)=0,i<j \right \}$, $1 \leq j \leq k$. Cook, Lyall and Magyar \cite{CookLyallMagyar}, see also the work written by Iosevich and Magyar \cite{IosevichMagyar}, suggest that due to our choice of the polynomials $P_{ij}$ we have $\df \omega_{F_j} (x_j) = C_{F_j} \df \sigma^{d-j}_{r_j}(x_j)$, $r_j>0$. Thus, we obtain the desired expression for $\delta(f)$.

\subsection{Step 3:} Given $\ve>0$ let $\psi_{\ve}(x)=\ve^{-d} \psi \pa{\frac{x}{\ve}} \geq 0$, where $\psi \geq 0$ is a smooth function such that its Fourier transform, $\h{\psi}$, is a  smooth compactly supported cut-off function, satisfying $\h{\psi}(0)=1$ and $0 \leq \h{\psi} \leq 1$. Let $\mu_{\ve} := \mu \ast \psi_{\ve}$. 

Here we will show that the continuous functions $\delta(\mu_{\ve})$ converge strongly (uniformly) as $\ve \to 0$. The rest of the proof of Theorem \ref{thm:maunthm1} relies on the following:
\begin{lem} \label{lemma1} There is a function $\dps M\pa{\f{t}} = \pa{\prod^{k}_{p=1} C_{F_p}}$ such that 
\benn
     \abs{\delta(\mu_{2\ve})-\delta(\mu_{\ve})} \lesssim |M \pa{\f{t}}| \ve^{\gamma}
\eenn
Where $\dps \gamma = \left\{\begin{matrix} (k-1)s - \frac{(k-1)(4k-1)d}{4k} - \frac{(k-1)}{4} & , & 3 \leq d < \frac{k(k+3)}{k-1} \\ (k-1)s - (k-1)d +1 & , &   \frac{k(k+3)}{k-1} \leq d \end{matrix}\right.$.
\end{lem}

%

\begin{rek}
We remind to the reader that $\dps C_{F_p}:= C_{F_p} \pa{\f{t}}$ are functions of $\f{t}$, more precisely, $\dps C_{F_p} = 2^{-p} \det (A_p)^{-1/2}$. Where $A_p = \pa{a_{ij}}$ is a $p \times p$ matrix, and $a_{ij}$ are described as in (\ref{eqn:mtdes}).  
\end{rek}

For the proof of Lemma \ref{lemma1} we build on the techniques developed by Iosevich and Magyar \cite{IosevichMagyar}. Details can be found in Section $4$. To have a better understanding of the proof of Lemma \ref{lemma1} the reader can go over Section $3$ in which we provide a proof for the case $k=2$.

Consider the sequence $\lla{\delta(\mu_{\ve_n})}_{n \in \N}$. Where $\ve_j = \frac{\ve}{2^j}$. Consider $m,n \in \N$, with $m>n$, then by Lemma \ref{lemma1} 
\begin{align*}
  \abs{\delta(\mu_{\ve_m}) - \delta(\mu_{\ve_n})} & \leq \sum^{m}_{j=n+1} \abs{\delta(\mu_{2\ve_j}) - \delta(\mu_{\ve_j})} \\
 & \lesssim \abs{M \pa{\f{t}}} \ve^{ \gamma } \sum^{m-1}_{j=n} 2^{-j \gamma  }\\
 & \lesssim \ve^{ \gamma } \sum^{\infty}_{j=1} 2^{-j \gamma }.
\end{align*} 
Note that the geometric series o the right-hand side is convergent as long as $\gamma > 0 $. From what we have
\benn
 \begin{matrix}
 s > \frac{(4k-1)d}{4k} + \frac{1}{4} & \text{if} & 3 \leq d < \frac{k(k+3)}{k-1} \\
 s > d - \frac{1}{k-1} & \text{if} & \frac{k(k+3)}{k-1} \leq d \\
\end{matrix}.
\eenn
Thus, the sequence $ \lla{\delta(\mu_{\ve_n})}_{n \in \N}$ is a Cauchy sequence with respect to the supremum norm. Therefore, $ \lla{\delta(\mu_{\ve_n})}_{n \in \N}$ converges uniformly to a continuous function, say $\delta(\mu)^{*}$. Due to the functions $\delta(\mu_{\ve})$ converge weakly to $\delta(\mu)$, then by uniqueness of the limit we have $\df \delta(\mu) = \delta(\mu)^{*} \df \f{t}$. Finally, we note that $\delta(\mu)$, as a pushforward of probability measures, is positive, and therefore the density $\delta(\mu)^{*}$ is non-zero. 

\qedd

\begin{rek}
    The reader might be concerned about the fact that support of $\mu_{\ve}$ might not be compact. However, $\mu_{\ve}$ is rapidly decreasing in a small neighborhood of the support of $\mu$, which is compact. This will not change the estimations given in Lemma \ref{lemma1} nor in the proof of our main result. For details see Remark \ref{remarkcs} at the end of Section $5$.
\end{rek}

\section{Proof of Lemma \ref{lemma1} for the case $k=2$.}
From (\ref{deltafs}) we have
\begin{multline*}
    \delta(\mu_{\ve})(t_{01},t_{02},t_{12}) = \iiint C_{F_1} C_{F_2} \mu_{\ve}(x) \mu_{\ve}(x+x_1)  \mu_{\ve}(x+x_2) \df \sigma^{d-2}_{r_2}(m x_1-x_2)  \df \sigma^{d-1}_{r_1}(x_1)  \df x.
\end{multline*} 
\begin{figure}[h!]
    \centering
    \includegraphics[scale=0.33]{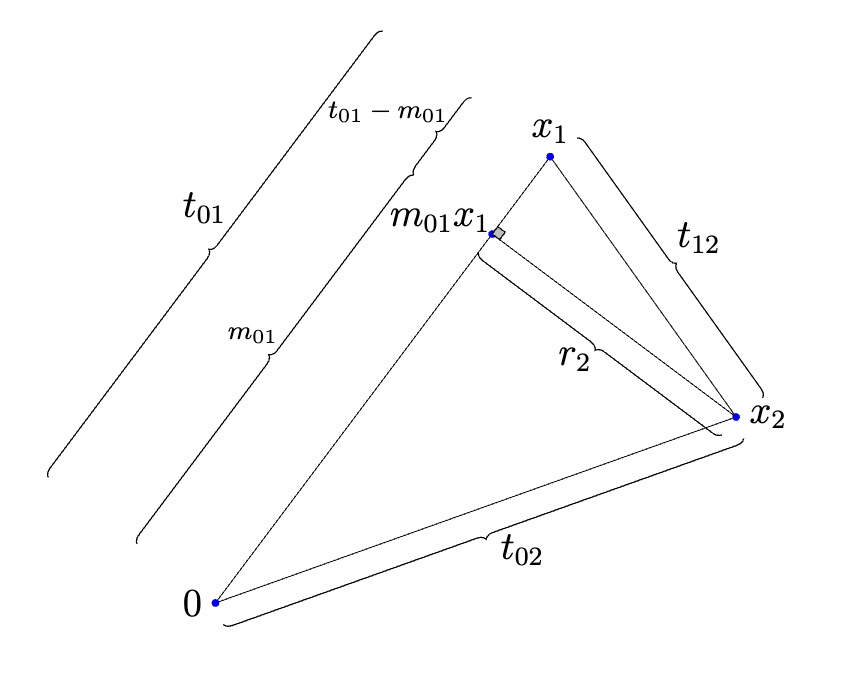}
	\caption{Here $C_{F_1}=\frac{t^2_{01}}{2}$, $C_{F_2}= \frac{1}{4} \corch{ t^2_{02} t^2_{12} - \frac{1}{4} \pa{t^2_{02} + t^2_{12} - t^2_{01}}^2 } $.}
	\label{fig:trian}
\end{figure}

\begin{figure}[h!]
    \centering
    \includegraphics[scale=0.33]{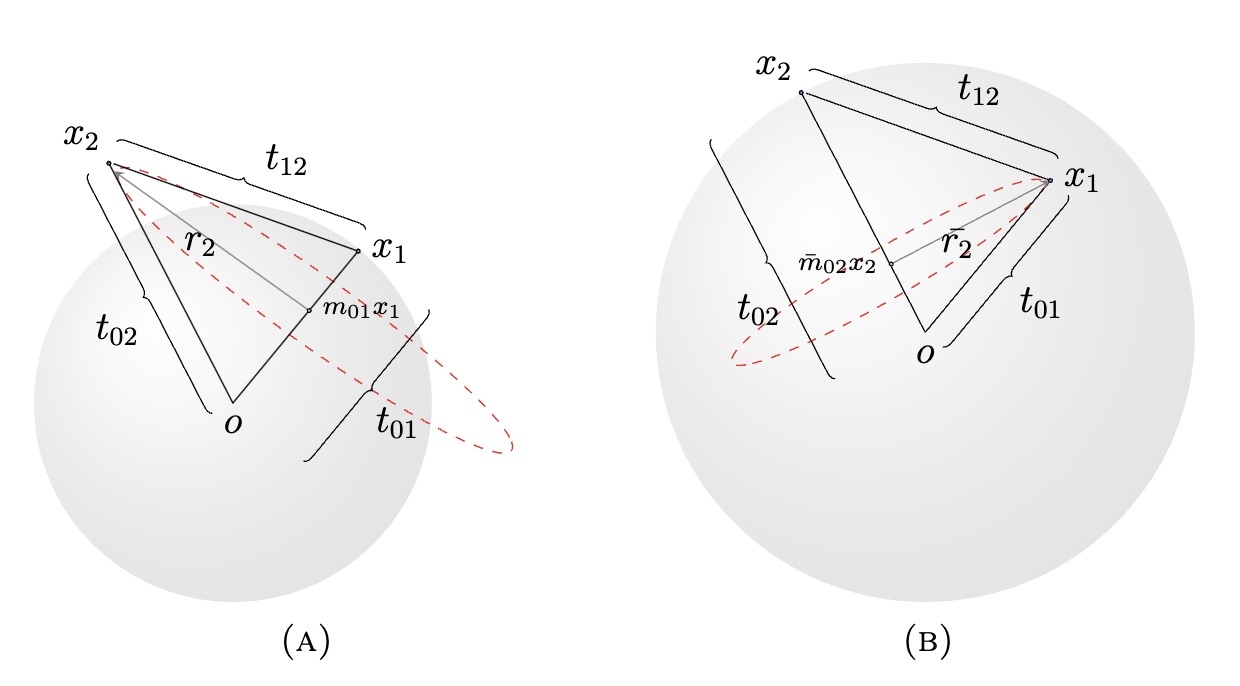}
	\caption{A triangle can be measured in different ways. In (A) we use the measure  $\df \sigma^{d-1}_{t_{01}}(x_1) \df \sigma^{d-2}_{r_2}(m_{01}x_1-x_2)$, but in (B) we use $\df \sigma^{d-1}_{t_{02}}(x_2) \df \sigma^{d-2}_{\Bar{r}_2}(\Bar{m}_{02}x_2-x_1)$. Thus, we can use either of these measures in the second term of the sum given in the right hand side of equation (\ref{triangleope}). A similar analysis can be done in each of the other terms.}
	\label{fig:esfs}
\end{figure}

Where $r_1=t_{01}$ , $m_{01}= \frac{t^2_{01}+t^2_{01}-t^2_{12}}{2t_{01}}$, and $r_2= t^2_{02} - m^2_{01}$ (see Figure \ref{fig:trian}). Write $\tri \mu_{\ve}:=\mu_{2\ve} - \mu_{\ve}$ , and $\tri \delta(\mu_{\ve}):=\delta(\mu_{2\ve}) - \delta(\mu_{\ve})$. Then
\benn 
    \prod^{2}_{j=0} \mu_{2\ve}(x+x_j) - \prod^{2}_{j=0} \mu_{\ve}(x+x_j) = \sum^{2}_{j=0} \tri_j(\mu_{\ve}), 
\eenn
where $\dps \tri_j(\mu_{\ve}) = \prod^{2}_{i=0,\ i \neq j } \mu_{\ve_{ij}}(x+x_i) \tri \mu_{\ve}(x+x_j) $, with $\ve_{ij} =  \left\{\begin{matrix}
2\ve &, i<j \\ \ve &, i>j \end{matrix}\right.$, and $x_0=0$. Therefore, 
\be\label{triangleope} 
\abs{\tri \delta(\mu_{\ve})} = \abs{C_{F_1} C_{F_2}} \abs{ \sum^{2}_{j=0}  \iiint \tri_j(\mu_{\ve}) \df \sigma^{d-2}_{r_2}(m_{01} x_1-x_2)  \df \sigma^{d-1}_{r_1}(x_1)  \df x }
\ee
 Note that by a simple change of variables or by considering the points in a different order, all the terms of the sum in the right-hand side of the equation (\ref{triangleope}) are equivalent (See Figure \ref{fig:esfs}). Therefore, we will just study the case $j=2$.
 
\begin{align*}
&  \abs{ \iiint \mu_{2\ve}(x) \mu_{2\ve}(x+x_1) \triangle \mu_{ \ve}(x+x_2) \df \sigma^{d-2}_{r_2}(x_2-m_{01} x_1) \df \sigma^{d-1}_{r_1}(x_1) \df x} \\
& \lesssim \ve^{s-d} \int \abs{ \iint \mu_{2\ve}(x) \triangle \mu_{ \ve}(x+x_2) \df \sigma^{d-2}_{r_2}(x_2-m_{01} x_1) \df x } \df \sigma^{d-1}_{r_1}(x_1)  \\
& = \ve^{s-d} \int \abs{ \int \mu_{2\ve}(x)  \corch{ \pa{\triangle \mu_{\ve} \ast \sigma^{d-2}_{r_2}} \circ \tau_{m_{01} x_1}}(x) \df x } \df \sigma^{d-1}_{r_1}(x_1)
\end{align*}
by Plancherel 
\begin{align*}
& \lesssim \ve^{s-d} \int \int \abs{ \h{\mu_{2\ve}}(\xi) }  \abs{ \h{ \triangle \mu_{\ve}}(\xi) } \abs{ \h{\sigma^{d-2}_{r_2}}(\xi) } \df \xi \df \sigma^{d-1}_{r_1}(x_1) \\
& = \ve^{s-d} \int \int \abs{ \h{\mu_{2\ve}}(\xi) } |\xi|^{\frac{1}{8}-\frac{d}{16}} |\xi|^{-\frac{1}{8}+\frac{d}{16}} \abs{ \h{ \triangle \mu_{\ve}}(\xi) } \abs{ \h{\sigma^{d-2}_{r_2}}(\xi) } \df \xi \df \sigma^{d-1}_{r_1}(x_1)
\end{align*}
and by applying Cauchy-Schwarz twice we have
\begin{align*}
& \leq \ve^{s-d} \int \pa{ \int \abs{ \h{\mu_{2\ve}}(\xi) }^2 |\xi|^{\frac{1}{4}-\frac{d}{8}} \df \xi }^{\frac{1}{2}} \pa{ \int |\xi|^{-\frac{1}{4}+\frac{d}{8}} \abs{ \h{ \triangle \mu_{\ve}}(\xi) }^2  \abs{ \h{\sigma^{d-2}_{r_2}}(\xi) }^2 \df \xi }^{\frac{1}{2}} \df \sigma^{d-1}_{r_1}(x_1) \\ 
& \lesssim  \ve^{s-d} \pa{ \int |\xi|^{-\frac{1}{4}+\frac{d}{8}} \abs{ \h{ \triangle \mu_{\ve}}(\xi)}^2 \pa{\int \abs{ \h{\sigma^{d-2}_{r_2} }(\xi) }^2 \df \sigma^{d-1}_{r_1}(x_1)} \df \xi }^{\frac{1}{2}}.
\end{align*}
Thus 
\benn 
  \abs{\triangle \delta(\mu_{\ve})}^2 \lesssim \abs{C_{F_1} C_{F_2}}^2 \ve^{2(s-d)} \pa{ \int |\xi|^{-\frac{1}{4}+\frac{d}{8}} \abs{ \h{ \triangle \mu_{\ve}}(\xi) }^2 T(\xi) \df \xi }. 
\eenn 
Where $\dps T(\xi) = \int \abs{ \h{\sigma^{d-2}_{r_2} }(\xi) }^2 \df \sigma^{d-1}_{r_1}(x_1)$. Due to $\sigma^{d-2}_{r_2}$ is $d-2$ dimensional sphere of radius $r_2>0$ is contained in an affine subspace orthogonal to the subspace $M=\lla{mx_1; \ m \in \R}$, then
\benn
 \abs{ \h{\sigma^{d-2}_{r_2} }(\xi)} \lesssim \pa{1+r_2 dist\pa{\xi,M}}^{-\frac{(d-2)}{2}}.
\eenn 
Note that the measure $\sigma^{d-1}_{r_1}(x_1)$ is invariant with respect to the change of variables $(x_1,x_2) \to (R(x_1),R(x_2))$ for any rotation $R \in SO(d)$, thus
\begin{align*}
  T(\xi) & \lesssim \iint \pa{1+r_2 dist\pa{\xi,M}}^{-\frac{(d-2)}{2}} \df \sigma^{d-1}_{r_1}(x_1) \df R \\
  & = \iint \pa{1+r_2 dist\pa{R(\xi),M}}^{-\frac{(d-2)}{2}} \df \sigma^{d-1}_{r_1}(x_1) \df R \\
  & = \iint \pa{1+r_2 |\xi|dist\pa{\eta,M}}^{-\frac{(d-2)}{2}} \df \sigma^{d-1}_{r_1}(x_1) \df \sigma^{d-2}(\eta) \\
  & \lesssim  \pa{1+r_2 |\xi|}^{-\frac{(d-2)}{2}}.
\end{align*}
Where $\eta := |\xi|^{-1}R(\xi)$. Then $\dps |T(\xi)| \lesssim |\xi|^{- \frac{(d-2)}{2}}$. Note that $\h{\tri \mu_{\ve}}(\xi)=\h{\mu}(\xi) \pa{ \h{\psi}(2\ve \xi) - \h{\psi}(\ve \xi) }$ is supported on $|\xi| \lesssim \ve^{-1}$. Thus,
\begin{align*}
 \int  |\xi|^{-\frac{1}{4}+\frac{d}{8}} \abs{ \h{ \triangle \mu_{\ve}}(\xi)}^2 T(\xi) \df \xi & = \int \limits_{|\xi| \leq \ve^{-1} }  |\xi|^{-\frac{1}{4}+\frac{d}{8}} \abs{ \h{ \triangle \mu_{\ve}}(\xi)}^2 T(\xi) \df \xi \\
    & = \sum^{ \floor{\log_2(\ve^{-1})} }_{j=-\infty} \int \limits_{  |\xi| \approx 2^j}  |\xi|^{-\frac{1}{4}+\frac{d}{8}} \abs{ \h{ \triangle \mu_{\ve}}(\xi)}^2 T(\xi) \df \xi.
\end{align*}
Observe that $\dps \abs{\h{\psi}(2\ve \xi) - \h{\psi}(\ve \xi)} \lesssim \ve |\xi| $. Thus the summation above is equal to
\begin{multline*}
  \sum^{ 0 }_{j=-\infty}  \ve^{2} \int \limits_{  |\xi| \approx 2^j}  |\xi|^{-\frac{1}{2}+\frac{d}{4}} |\xi|^{\frac{1}{4}-\frac{d}{8}} \abs{ \h{\mu}(\xi)}^2 |\xi|^2 T(\xi) \df \xi  \\ 
  + \sum^{ \floor{\log_2(\ve^{-1})} }_{j= 1} \ve^{2} \int \limits_{  |\xi| \approx 2^j}   |\xi|^{-\frac{1}{2}+\frac{d}{4}} |\xi|^{\frac{1}{4}-\frac{d}{8}} \abs{ \h{\mu}(\xi)}^2 |\xi|^2 T(\xi) \df \xi.   
\end{multline*}
Note
\begin{multline*}
    \sum^{ 0 }_{j=-\infty}  \ve^{2} \int \limits_{  |\xi| \approx 2^j}  |\xi|^{-\frac{1}{2}+\frac{d}{4}} |\xi|^{\frac{1}{4}-\frac{d}{8}} \abs{ \h{\mu}(\xi)}^2 |\xi|^2 T(\xi) \df \xi \\
    \lesssim \ve^2 \sum^{0}_{j = - \infty } (2^j)^{-\frac{1}{2}+\frac{d}{4}+2} \int \limits_{  |\xi| \approx 2^j}   |\xi|^{\frac{1}{4}-\frac{d}{8}} \abs{\h{ \mu}(\xi)}^2 \df \xi \lesssim \ve^2
\end{multline*}
and
\begin{multline} \label{eqn:bound1}
 \sum^{  \floor{\log_2(\ve^{-1})} }_{j=1}  \ve^{2} \int \limits_{  |\xi| \approx 2^j}  |\xi|^{-\frac{1}{2}+\frac{d}{4}} |\xi|^{\frac{1}{4}-\frac{d}{8}} \abs{ \h{\mu}(\xi)}^2 |\xi|^2 T(\xi) \df \xi 
 \\ \lesssim  \ve^2 \sum^{\floor{\log_2(\ve^{-1})}}_{j = 1 } (2^j)^{-\frac{1}{2}+\frac{d}{4}+2 -\frac{(d-2)}{2} } \int \limits_{  |\xi| \approx 2^j} |\xi|^{\frac{1}{4}-\frac{d}{8}} \abs{\h{ \mu}(\xi)}^2 \df \xi 
\end{multline} 
Observe that there are two ways in which we can bound the right-hand side of inequality (\ref{eqn:bound1}). We can add the terms to infinity, that is, 
\begin{multline*}
 \sum^{  \floor{\log_2(\ve^{-1})} }_{j=1}  \ve^{2} \int \limits_{  |\xi| \approx 2^j}  |\xi|^{-\frac{1}{2}+\frac{d}{4}} |\xi|^{\frac{1}{4}-\frac{d}{8}} \abs{ \h{\mu}(\xi)}^2 |\xi|^2 T(\xi) \df \xi 
 \\ \lesssim  \ve^2 \sum^{ \infty }_{j = 1 } \pa{2^{-\frac{1}{2}+\frac{d}{4}+2 -\frac{(d-2)}{2} }}^j \int \limits_{ |\xi| \approx 2^j} |\xi|^{\frac{1}{4}-\frac{d}{8}} \abs{\h{ \mu}(\xi)}^2 \df \xi  \lesssim \ve^2,
\end{multline*} 
which requires $d>10$. We can also add the terms of the finite sum, that is, 
\begin{align*}
 & \sum^{  \floor{ \log_2(\ve^{-1}) } }_{j=1}  \ve^{2} \int \limits_{  |\xi| \approx 2^j}  |\xi|^{-\frac{1}{2}+\frac{d}{4}} |\xi|^{\frac{1}{4}-\frac{d}{8}} \abs{ \h{\mu}(\xi)}^2 |\xi|^2 T(\xi) \df \xi \\
 & \lesssim  \ve^2 \sum^{ \floor{\log_2(\ve^{-1})}  }_{j = 1 } \pa{2^{-\frac{1}{2}+\frac{d}{4}+2 -\frac{(d-2)}{2} }}^j \int \limits_{ |\xi| \approx 2^j} |\xi|^{\frac{1}{4}-\frac{d}{8}} \abs{\h{ \mu}(\xi)}^2 \df \xi \\
 &\lesssim \ve^2 \corch{ \frac{  \ve^{ -\pa{ \frac{5}{2} -\frac{d}{4}} } - 1 }{\pa{2^{ \frac{5}{2} -\frac{d}{4} } }  -1 } }.
\end{align*} 
In the first case we then have $\dps \abs{\triangle \delta(\mu_{\ve})}^2 \lesssim \abs{C_{F_1} C_{F_2}}^2 \ve^{2(s-d)} \pa{ 2 \ve^2 }$, and thus
\benn 
 \abs{\triangle \delta(\mu_{\ve})} \lesssim \abs{C_{F_1} C_{F_2}} \ve^{s-d+1 }, 
\eenn
if $d>10$. In the second case we have 
\benn
\abs{\triangle \delta(\mu_{\ve})}^2 \lesssim \abs{C_{F_1} C_{F_2}}^2 \ve^{2(s-d)} \pa{ \ve^2 + \ve^2 \corch{ \frac{  \ve^{ -\pa{ \frac{5}{2} -\frac{d}{4}} } - 1 }{\pa{2^{ \frac{5}{2} -\frac{d}{4} } }  -1 } } }
\eenn 
thus
\begin{align*}
    \abs{\triangle \delta(\mu_{\ve})} & \lesssim \abs{C_{F_1} C_{F_2}} \ve^{ \frac{1}{2} \pa{ 2(s-d)+2 -\pa{ \frac{5}{2} -\frac{d}{4}} }  } \pa{ \ve^{ \frac{5}{2} -\frac{d}{4} }  +   \corch{ \frac{ 1 - \ve^{\frac{5}{2} -\frac{d}{4}} }{\pa{2^{ \frac{5}{2} -\frac{d}{4} } } -1 } } }^{ \frac{1}{2} } \\
    & \lesssim \abs{C_{F_1} C_{F_2}} \ve^{ s - \frac{7d}{8} - \frac{1}{4} },   
\end{align*}
if $d<10$. When $d=10$, we use the following  equation
\begin{multline*}
     \int  |\xi|^{-\frac{1}{4}+\frac{10}{8}} \abs{ \h{ \triangle \mu_{\ve}}(\xi)}^2 T(\xi) \df \xi = \int \limits_{|\xi| \leq \ve  }  |\xi|^{-\frac{1}{4}+\frac{10}{8}} \abs{ \h{ \triangle \mu_{\ve}}(\xi)}^2 T(\xi) \df \xi \\ + \int \limits_{\ve  < |\xi| \leq \ve^{-1} }  |\xi|^{-\frac{1}{4}+\frac{10}{8}} \abs{ \h{ \triangle \mu_{\ve}}(\xi)}^2 T(\xi) \df \xi
\end{multline*}
from which we have 
\benn 
    \int  |\xi|^{-\frac{1}{4}+\frac{10}{8}} \abs{ \h{ \triangle \mu_{\ve}}(\xi)}^2 T(\xi) \df \xi \lesssim \ve^2 \pa{\ve^4+1 },
\eenn
thus $\dps \abs{\triangle \delta(\mu_{\ve})}^2 \lesssim \abs{C_{F_1} C_{F_2}}^2 \ve^{2(s-10)} \ve^2 \pa{  \ve^4 +1 }$, and therefore
\benn 
 \abs{\triangle \delta(\mu_{\ve})} \lesssim \abs{C_{F_1} C_{F_2}} \ve^{s-9 }. 
\eenn
In other words we have 
\benn 
    \abs{\triangle \delta(\mu_{\ve})} \lesssim \abs{C_{F_1} C_{F_2}}  \ve^{\gamma}.
\eenn
Where $\dps \gamma = \left\{\begin{matrix} s - \frac{7d}{8} - \frac{1}{4} & ,\text{if} & 3 \leq d < 10 \\ s-d+1 & ,\text{if} &   10 \leq d \end{matrix}\right.$.

\qedd

\section{Proof of Lemma \ref{lemma1} general case}
In this section we will show that the proof given above can be easily extended. From (\ref{deltafs}) we have
\begin{multline*}
 \delta(\mu_{\ve})\pa{\f{t}} = \pa{\prod^{k}_{p=1} C_{F_p}} \int \hdots \int \mu_{\ve}(x) \mu_{\ve}(x+x_1) \dots \mu_{\ve}(x+x_k) \\ \df \sigma^{d-k}_{r_k}\pa{ x_k-\sum^{k-1}_{n=1} m_{nk}x_n} \hdots \df \sigma^{d-2}_{r_2}(x_2-m_{12}x_1)  \df \sigma^{d-1}_{r_1}(x_1) \df x.
\end{multline*} 
Write $\tri \mu_{\ve}:=\mu_{2\ve} - \mu_{\ve}$ , and $\tri \delta(\mu_{\ve}):=\delta(\mu_{2\ve}) - \delta(\mu_{\ve})$. Then
\benn 
    \prod^{k}_{j=0} \mu_{2\ve}(x+x_j) - \prod^{k}_{j=0} \mu_{\ve}(x+x_j) = \sum^{k}_{j=0} \tri_j(\mu_{\ve})
\eenn
where $\dps \tri_j(\mu_{\ve}) = \prod^{k}_{i=0,\ i \neq j } \mu_{\ve_{ij}}(x+x_i) \tri \mu_{\ve}(x+x_j) $, with $\dps \ve_{ij} = \left\{ \begin{matrix} 2\ve &, i<j \\ \ve &, i>j \end{matrix} \right. $, and $x_0=0$. Therefore, 
\be \label{triangleopes} 
\abs{\tri \delta(\mu_{\ve})} = \abs{\prod^{k}_{p=1} C_{F_p}} \abs{ \sum^{k}_{j=0} \int \hdots \int \tri_j(\mu_{\ve}) \prod^{k}_{i=1} \df \sigma^{d-i}_{r_i}\pa{ x_i-\sum^{i-1}_{n=0} m_{ni}x_n} \df x }
\ee
By a similar reasoning as in the proof of Lemma \ref{lemma1} for the case $k=2$, we can conclude that all the terms of the sum in the right hand side of equation (\ref{triangleopes}) are equivalent, therefore we just study the case $j=k$. Note
 \begin{multline*}
 \left | \int \hdots \int \prod^{k-1}_{i=0} \mu_{\ve_{ij}}(x+x_i) \tri \mu_{\ve}(x+x_k) \df \sigma^{d-k}_{r_k}(x_k-m x_1) \right. \\ \left.  \prod^{k}_{i=1} \df \sigma^{d-i}_{r_i}\pa{ x_i-\sum^{i-1}_{n=0} m_{ni}x_n} \df x \right |
 \end{multline*}
  \begin{multline*}
  \lesssim  \ve^{(k-1)(s-d)} \int \abs{ \iint \mu_{2\ve}(x) \triangle \mu_{ \ve}(x+x_k) \df \sigma^{d-k}_{r_k}\pa{ x_k-\sum^{k-1}_{n=0} m_{nk}x_n} \df x } \\ \df \Bar{\omega} \pa{x_1,x_2, \hdots,x_{k-1}}
  \end{multline*}
  \begin{multline*}
  = \ve^{(k-1)(s-d)} \int \abs{ \int \mu_{2\ve}(x) \corch{ \pa{\triangle \mu_{\ve} \ast \sigma^{d-k}_{r_k}} \circ \tau_{\Bar{x}}}(x) \df x } \df \Bar{\omega} \pa{x_1,x_2, \hdots,x_{k-1}}.
  \end{multline*}
Where $\Bar{x}=\sum^{k-1}_{n=0} m_{nk}x_n$, and $\dps \df \Bar{\omega} \pa{x_1,x_2, \hdots,x_{k-1}} = \prod^{k-1}_{i=1} \df \sigma^{d-i}_{r_i}\pa{ x_i-\sum^{i-1}_{n=0} m_{ni}x_n}$. By Plancherel
\begin{align*}
 & \ve^{(k-1)(s-d)} \int \abs{ \int \mu_{2\ve}(x) \corch{ \pa{\triangle \mu_{\ve} \ast \sigma^{d-k}_{r_k}} \circ \tau_{\Bar{x}}}(x) \df x } \df \Bar{\omega} \pa{x_1,x_2, \hdots,x_{k-1}}\\
 & \leq  \ve^{(k-1)(s-d)} \int \int \abs{ \h{\mu_{2\ve}}(\xi) }  \abs{ \h{ \triangle \mu_{\ve}}(\xi) }  \abs{ \h{\sigma^{d-k}_{r_k}}(\xi) } \df \xi \df \Bar{\omega} \pa{x_1,x_2, \hdots,x_{k-1}} \\
  & =  \ve^{(k-1)(s-d)} \int \int \abs{ \h{\mu_{2\ve}}(\xi) }  |\xi|^{\frac{1}{8}-\frac{d}{8k}} |\xi|^{-\frac{1}{8}+\frac{d}{8k}} \abs{ \h{ \triangle \mu_{\ve}}(\xi) }  \abs{ \h{\sigma^{d-k}_{r_k}}(\xi) } \df \xi \df \Bar{\omega} \pa{x_1,x_2, \hdots,x_{k-1}} 
\end{align*}
and by applying Cauchy-Schwarz twice we have
\begin{multline*}
\leq \ve^{(k-1)(s-d)} \int \pa{ \int \abs{ \h{\mu_{2\ve}}(\xi) }^2 |\xi|^{\frac{1}{4}-\frac{d}{4k}} \df \xi }^{\frac{1}{2}} \pa{ \int |\xi|^{-\frac{1}{4}+\frac{d}{4k}} \abs{ \h{ \triangle \mu_{\ve}}(\xi) }^2  \abs{ \h{\sigma^{d-k}_{r_k}}(\xi) }^2 \df \xi }^{\frac{1}{2}} \\ 
\df \Bar{\omega} \pa{x_1,x_2, \hdots,x_{k-1}} 
\end{multline*}
\benn
 \lesssim \ve^{(k-1)(s-d)} \pa{ \int |\xi|^{-\frac{1}{4}+\frac{d}{4k}} \abs{ \h{ \triangle \mu_{\ve}}(\xi) }^2 \pa{ \int \abs{ \h{\sigma^{d-k}_{r_k} }(\xi) }^2 \df \Bar{\omega} \pa{x_1,x_2, \hdots,x_{k-1}} } \df \xi }^{\frac{1}{2}}  
\eenn
Thus 
\benn
  \abs{\triangle \delta(\mu_{\ve})}^2 \lesssim \abs{\prod^{k}_{p=1} C_{F_p}}^2 \ve^{2(k-1)(s-d)}  \pa{ \int |\xi|^{-\frac{1}{4}+\frac{d}{4k}} \abs{ \h{ \triangle \mu_{\ve}}(\xi) }^2 T(\xi) \df \xi }. 
 \eenn
Where $\dps T(\xi) = \int \abs{ \h{\sigma^{d-k}_{r_k} }(\xi) }^2 \df \Bar{\omega} \pa{x_1,x_2, \hdots,x_{k-1}} $. Let $r=\min \lla{r_n, \ 1 \leq n \leq k}$. Due to $\sigma^{d-k}_{r_k} $ is a $d-k$ dimensional sphere of radius $r_k \geq r>0$ is contained in an affine subspace orthogonal to the subspace $M(x_1,\hdots,x_{k-1})= Span \lla{x_1,x_2, \hdots ,x_{k-1} }$, then
\benn
 \abs{ \h{\sigma^{d-k}_{r_k}}(\xi)} \lesssim \pa{1+rdist\pa{\xi,M(x_1,\hdots,x_{k-1})}}^{-\frac{(d-k)}{2}}.
\eenn 
Note that the measure $\df \Bar{\omega} \pa{x_1,x_2, \hdots,x_{k-1}}$ is invariant with respect to the change of variables $(x_1,\hdots, x_{k-1}) \to (R(x_1),\hdots,R(x_{k-1}))$ for any rotation $R \in SO(d)$, thus
\begin{align*}
  T(\xi) & \lesssim  \iint \pa{1+r dist\pa{\xi,M(R(x_1),\hdots,R(x_{k-1}))}}^{-\frac{(d-k)}{2}} \df \Bar{\omega} \pa{x_1,x_2, \hdots,x_{k-1}} \df R \\
  & =  \iint \pa{1+r dist\pa{R(\xi),M(x_1,\hdots,x_{k-1})}}^{-\frac{(d-k)}{2}} \df \Bar{\omega} \pa{x_1,x_2, \hdots,x_{k-1}} \df R \\
  & =  \iint \pa{1+r |\xi| dist\pa{\eta,M(x_1,\hdots,x_{k-1})}}^{-\frac{(d-k)}{2}} \df \Bar{\omega} \pa{x_1,x_2, \hdots,x_{k-1}} \df \sigma^{k-1}(\eta)\\
  & \lesssim  \pa{1+r|\xi|}^{-\frac{(d-k)}{2}}.
\end{align*}
Where $\eta := |\xi|^{-1}R(\xi)$. Then $\dps |T(\xi)| \lesssim |\xi|^{- \frac{(d-k)}{2}}$. Note that $\h{\tri \mu_{\ve}}(\xi)=\h{\mu}(\xi) \pa{ \h{\psi}(2\ve \xi) - \h{\psi}(\ve \xi) }$ is supported on $|\xi| \lesssim \ve^{-1}$. Thus,
\begin{align*}
 \int  |\xi|^{-\frac{1}{4}+\frac{d}{4k}} \abs{ \h{ \triangle \mu_{\ve}}(\xi)}^2 T(\xi) \df \xi & = \int \limits_{|\xi| \leq \ve^{-1} }  |\xi|^{-\frac{1}{4}+\frac{d}{4k}} \abs{ \h{ \triangle \mu_{\ve}}(\xi)}^2 T(\xi) \df \xi \\
    & = \sum^{ \floor{\log_2(\ve^{-1})} }_{j=-\infty} \int \limits_{  |\xi| \approx 2^j}  |\xi|^{-\frac{1}{4}+\frac{d}{4k}} \abs{ \h{ \triangle \mu_{\ve}}(\xi)}^2 T(\xi) \df \xi.
\end{align*}
Observe that $\dps \abs{\h{\psi}(2\ve \xi) - \h{\psi}(\ve \xi)} \lesssim \ve |\xi| $. Thus the summation above is equal to
\begin{multline*}
  \sum^{ 0 }_{j=-\infty}  \ve^{2} \int \limits_{  |\xi| \approx 2^j}  |\xi|^{-\frac{1}{2}+\frac{d}{2k}}  |\xi|^{\frac{1}{4}-\frac{d}{4k}} \abs{ \h{\mu}(\xi)}^2 |\xi|^2 T(\xi) \df \xi  \\ 
  + \sum^{ \floor{\log_2(\ve^{-1})} }_{j= 1} \ve^{2} \int \limits_{  |\xi| \approx 2^j}   |\xi|^{-\frac{1}{2}+\frac{d}{2k}}  |\xi|^{\frac{1}{4}-\frac{d}{4k}} \abs{ \h{\mu}(\xi)}^2 |\xi|^2 T(\xi) \df \xi
\end{multline*}
Note
\begin{multline*}
    \sum^{ 0 }_{j=-\infty}  \ve^{2} \int \limits_{  |\xi| \approx 2^j}  |\xi|^{-\frac{1}{2}+\frac{d}{2k}} |\xi|^{\frac{1}{4}-\frac{d}{4k}} \abs{ \h{\mu}(\xi)}^2 |\xi|^2 T(\xi) \df \xi \\
    \lesssim \ve^2 \sum^{0}_{j = - \infty } (2^j)^{-\frac{1}{2}+\frac{d}{2k}+2} \int \limits_{  |\xi| \approx 2^j}   |\xi|^{\frac{1}{4}-\frac{d}{4k}} \abs{\h{ \mu}(\xi)}^2 \df \xi \lesssim \ve^2
\end{multline*}
and
\begin{multline} \label{eqn:bound2}
 \sum^{  \floor{\log_2(\ve^{-1})} }_{j=1}  \ve^{2} \int \limits_{  |\xi| \approx 2^j}  |\xi|^{-\frac{1}{2}+\frac{d}{2k}}  |\xi|^{\frac{1}{4}-\frac{d}{4k}} \abs{ \h{\mu}(\xi)}^2 |\xi|^2 T(\xi) \df \xi 
 \\ \lesssim  \ve^2 \sum^{\floor{\log_2(\ve^{-1})}}_{j = 1 } (2^j)^{-\frac{1}{2}+\frac{d}{2k} +2 -\frac{(d-k)}{2} } \int \limits_{  |\xi| \approx 2^j} |\xi|^{\frac{1}{4}-\frac{d}{4k}} \abs{\h{ \mu}(\xi)}^2 \df \xi 
\end{multline} 
Similar to the case $k=2$, there are two ways in which we can bound the right-hand side of inequality (\ref{eqn:bound2}). We can add the terms to infinity, that is, 
\begin{multline*}
 \sum^{  \floor{\log_2(\ve^{-1})} }_{j=1}  \ve^{2} \int \limits_{  |\xi| \approx 2^j}  |\xi|^{-\frac{1}{2}+\frac{d}{2k}}  |\xi|^{\frac{1}{4}-\frac{d}{4k}} \abs{ \h{\mu}(\xi)}^2 |\xi|^2 T(\xi) \df \xi 
 \\ \lesssim  \ve^2 \sum^{ \infty }_{j = 1 } \pa{2^{-\frac{1}{2}+\frac{d}{2k}+2 -\frac{(d-2)}{2} }}^j \int \limits_{ |\xi| \approx 2^j} |\xi|^{\frac{1}{4}-\frac{d}{4k}} \abs{\h{ \mu}(\xi)}^2 \df \xi  \lesssim \ve^2,
\end{multline*} 
which requires $d> \frac{k(k+3)}{k-1}$. We can also add the terms of the finite sum, that is, 
\begin{align*}
 & \sum^{  \floor{ \log_2(\ve^{-1}) } }_{j=1}  \ve^{2} \int \limits_{  |\xi| \approx 2^j}  |\xi|^{-\frac{1}{2}+\frac{d}{2k}}  |\xi|^{\frac{1}{4}-\frac{d}{4k}} \abs{ \h{\mu}(\xi)}^2 |\xi|^2 T(\xi) \df \xi \\
 & \lesssim  \ve^2 \sum^{ \floor{\log_2(\ve^{-1})}  }_{j = 1 } \pa{2^{-\frac{1}{2}+\frac{d}{2k}+2 -\frac{(d-k)}{2} }}^j \int \limits_{ |\xi| \approx 2^j} |\xi|^{\frac{1}{4}-\frac{d}{4k}} \abs{\h{ \mu}(\xi)}^2 \df \xi \\
 &\lesssim \ve^2 \corch{ \frac{  \ve^{ -\pa{ \frac{(1-k)d}{2k} + \frac{k+3}{2}} } - 1 }{\pa{2^{ \frac{(1-k)d}{2k} + \frac{k+3}{2} } }  -1 } }.
\end{align*} 
In the first case we then have $\dps \abs{\triangle \delta(\mu_{\ve})}^2 \lesssim \abs{\prod^{k}_{p=1} C_{F_p}}^2 \ve^{2(k-1)(s-d)} \pa{ 2 \ve^2 }$, and thus
\benn 
 \abs{\triangle \delta(\mu_{\ve})} \lesssim \abs{\prod^{k}_{p=1} C_{F_p}} \ve^{(k-1)s-(k-1)d+1 }, 
\eenn
if $d>\frac{k(k+3)}{k-1}$. In the second case we have 
\benn
\abs{\triangle \delta(\mu_{\ve})}^2 \lesssim \abs{\prod^{k}_{p=1} C_{F_p}}^2 \ve^{2(k-1)(s-d)} \pa{ \ve^2 + \ve^2 \corch{ \frac{  \ve^{ -\pa{ \frac{(1-k)d}{2k} + \frac{k+3}{2} } } - 1 }{\pa{2^{ \frac{(1-k)d}{2k} + \frac{k+3}{2} } }  -1 } } }
\eenn 
thus
\begin{align*}
    \abs{\triangle \delta(\mu_{\ve})} & \lesssim \abs{\prod^{k}_{p=1} C_{F_p}} \ve^{ \frac{1}{2} \pa{ 2(k-1)(s-d)+2 -\pa{ \frac{(1-k)d}{2k} + \frac{k+3}{2} } }  } \pa{ \ve^{ \frac{(1-k)d}{2k} + \frac{k+3}{2} }  +   \corch{ \frac{ 1 - \ve^{ \frac{(1-k)d}{2k} + \frac{k+3}{2} } }{\pa{2^{ \frac{(1-k)d}{2k} + \frac{k+3}{2} } } -1 } } }^{ \frac{1}{2} } \\
    & \lesssim \abs{\prod^{k}_{p=1} C_{F_p}} \ve^{ (k-1)s - \frac{(k-1)(4k-1)d}{4k} - \frac{(k-1)}{4} },   
\end{align*}
if $d<\frac{k(k+3)}{k-1}$. When $d=\frac{k(k+3)}{k-1}$, we use the following  equation
\begin{multline*}
     \int  |\xi|^{-\frac{1}{4}+\frac{1}{4k} \pa{ \frac{k(k+3)}{k-1}} } \abs{ \h{ \triangle \mu_{\ve}}(\xi)}^2 T(\xi) \df \xi = \int \limits_{|\xi| \leq \ve  }  |\xi|^{-\frac{1}{4}+\frac{1}{4k}  \pa{ \frac{k(k+3)}{k-1}} } \abs{ \h{ \triangle \mu_{\ve}}(\xi)}^2 T(\xi) \df \xi \\ + \int \limits_{\ve  < |\xi| \leq \ve^{-1} }  |\xi|^{-\frac{1}{4}+\frac{1}{4k}  \pa{ \frac{k(k+3)}{k-1}} } \abs{ \h{ \triangle \mu_{\ve}}(\xi)}^2 T(\xi) \df \xi
\end{multline*}
from which we have 
\benn 
    \int  |\xi|^{-\frac{1}{4}+\frac{1}{4k} \pa{ \frac{k(k+3)}{k-1}} } \abs{ \h{ \triangle \mu_{\ve}}(\xi)}^2 T(\xi) \df \xi \lesssim \ve^2 \pa{ \ve^{2 + \frac{2}{k-1}} +1 },
\eenn
thus $\dps \abs{\triangle \delta(\mu_{\ve})}^2 \lesssim \abs{\prod^{k}_{p=1} C_{F_p}}^2 \ve^{2(k-1)\pa{ s- \pa{ \frac{k(k+3)}{k-1}} } } \ve^2 \pa{  \ve^{2 + \frac{2}{k-1}} +1 }$, and therefore
\benn 
 \abs{\triangle \delta(\mu_{\ve})} \lesssim \abs{\prod^{k}_{p=1} C_{F_p}} \ve^{(k-1)s-k(k+3)+1 }. 
\eenn
In other words we have 
\benn 
    \abs{\triangle \delta(\mu_{\ve})} \lesssim \abs{\prod^{k}_{p=1} C_{F_p}} \ve^{\gamma}.
\eenn
Where $\dps \gamma = \left\{\begin{matrix} (k-1)s - \frac{(k-1)(4k-1)d}{4k} - \frac{(k-1)}{4} & , \text{if} & 3 \leq d < \frac{k(k+3)}{k-1} \\ (k-1)s - (k-1)d +1 & , \text{if} &  \frac{k(k+3)}{k-1} \leq d \end{matrix}\right.$.

\qedd

\begin{rek} \label{remarkcs}
    Although the support of $\mu_{\ve}$ might not be compact, it is rapidly decreasing in a small neighborhood of the support of $\mu$. Consider $\phi_{\ve}(x) := \phi \pa{ \ve^{-1/2}x}$, where $0 \leq \phi(x) \leq 1$ is a smooth cut-off function, such that 
     $$ \phi(x) = \left\{\begin{matrix} 1 &, \text{if} & |x| \leq \frac{1}{2} \\ 0 &, \text{if} & |x| \geq 2  \end{matrix}\right. ,$$
    and consider $ \widetilde{\psi}_{\ve} = \psi_{\ve} \phi_{\ve}$. Let $\widetilde{\mu}_{\ve} = \mu \ast \widetilde{\psi}_{\ve}$, and note that $\widetilde{\mu}_{\ve} \leq \mu_{\ve}$. Thus,
    \begin{align*}
        \abs{\delta(\mu) - \delta(\widetilde{\mu}_{\ve})} & \leq C_k M\pa{\f{t}} \nor{ \mu_{\ve}}^{k-1}_{\infty} \pa{ \int \abs{ \mu_{\ve} - \widetilde{\mu}_{\ve} } \df \Bar{\omega} \pa{x_1,x_2, \hdots,x_{k-1}} } \\
        & \leq C_k M\pa{\f{t}} \nor{ \mu_{\ve}}^{k-1}_{\infty} \nor{  \mu_{\ve} - \widetilde{\mu}_{\ve} }_{\infty} \\
         & \leq C_k M\pa{\f{t}} \ve^{ (s-d)(k-1) + \frac{N}{2} - d}.
    \end{align*}
    Where the last inequality comes from  
    \begin{align*} 
    \abs{ \mu_{\ve} - \widetilde{\mu}_{\ve} } & \leq \int \abs{ \ve^{-d} \psi\pa{\frac{y-x}{\ve}} } \abs{ 1- \phi \pa{\frac{y-x}{\ve^{1/2}}} } \df \mu(x) 
    \end{align*}
    Since $\psi$ is a Schwartz function we have $\dps \abs{ \psi\pa{\frac{y-x}{\ve}} } \lesssim_N \abs{  \frac{y-x}{\ve} }^{-N}$, then
    \begin{align*} 
    \abs{ \mu_{\ve} - \widetilde{\mu}_{\ve} } & \lesssim_{N} \int \ve^{-d} \abs{  \frac{y-x}{\ve} }^{-N}  \abs{ 1- \phi \pa{\frac{y-x}{\ve^{1/2}}} } \df \mu(x) \\
    & \lesssim_{N} \int_{ \ve^{1/2} \lesssim |y-x| } \ve^{-d} \abs{  \frac{y-x}{\ve} }^{-N} \df \mu(x) \\
    & \lesssim_{N} \ve^{ \frac{N}{2} - d} 
    \end{align*}
    thus, by taking $N$ large enough, we have
    $\dps \abs{\delta(\mu) - \delta(\widetilde{\mu}_{\ve})} \lesssim \ve$. The reader can easily show that the estimates given in the proof of Lemma \ref{lemma1} still hold if one replaces $\mu_{\ve}$ by $\widetilde{\mu}_{\ve}$. 
\end{rek}


\section{Proof of Lemma \ref{PHprinc2} }
Similar to the proof of Lemma \ref{PHprinc} given in \cite{PalssonRomero} (and Lemma $2.1$ in \cite{IosevichMourgoglouSenger}) we will use a stopping time argument to show that it is possible to find at least $d+1$ cubes with positive measure. 

Without lost of generality assume that $E \subset \corch{0,1}^{d}$ such that $\mu(E)=1$. Where $\corch{0,1}^{d}$ is the unit cube in $\R^d$ and $\mu$ is a Frostman probability measure. Let $C_{\mu}>0$ be the constant in the Frostman condition $\mu \pa{B(x,r)} \leq C_{\mu} r^s$, for some $s>0$ and for all $x \in \R^d$ and $r>0$. For the proof of this Lemma we may assume that $C_{\mu}$ is big enough.

Lets divide the unit cube $\corch{0,1}^{d}$ into $6^d$ smaller cubes with edge-length $\frac{1}{6}$. Let $\Omega_k$, $1 \leq k \leq 3^{d-1}2^2$, be a collection of $32^{d-2}$ sub-cubes such that no two cubes of the same collection touch each other. By pigeon hole principle at least one of the collections $\Omega_k$ has measure greater or equal to $ \frac{1}{3^{d-1}2^2}$, that is, 
\begin{center}
  $\dps \mu \left ( \bigcup_{Q \in \Omega_k} Q  \right ) \geq \frac{1}{3^{d-1}2^2}$ for some $k$.
\end{center}
We have the following cases:
\begin{itemize}
    \item[(1)] There is a collection that contains $d+1$ cubes with positive measure. If this is the case then we conclude the proof. 
    \item[(2)] For each $1 \leq l \leq d$, there is a collection that contains $l$ cubes with positive measure. In other words, for some $1\leq k \leq 3^{d-1}2^2$ there are cubes $Q_{ik} \in \Omega_k$ such that $\mu(Q_{ik}) > \frac{c}{3^{d-1}2^2}$, for $1 \leq i \leq l$ and some $c>0$. Then we have the following sub-cases
     \begin{itemize}
     \item [(a)] If $c \geq 1$, then we have $\mu (Q_{ik}) \geq \frac{1}{3^{d-1}2^2}$, for $1 \leq i \leq l$. Thus we pick one of these cubes, say $Q_{1k}$, and we repeat the procedure, that is, we subdivide $Q_{1k}$ into $3^{d-1}2^2$ collections of $32^{d-2}$ smaller cubes each.  
     \item [(b)] If $c<1$, then there is a cube, say $Q_{1k}$, such that $\mu (Q_{1k}) \geq \frac{1}{l3^{d-1}2^2}$, and thus we repeat the procedure on $Q_{1k}$. If such cube does not exists, then $0 < \mu (Q_{ik}) < \frac{1}{l3^{d-1}2^2}$ for all $1 \leq i \leq l$, therefore
            \begin{center}
              $\dps \mu \pa{ \bigcup_{Q \in \Omega_k \setminus \left \{ Q_{ik} \right \} } Q } >0 $.
            \end{center}
            Thus, there must be a cube in $\Omega_k \setminus \left \{ Q_{ik} \right \}$ with positive $\mu$ measure, and we have $l+1$ cubes with positive $\mu$ measure.  
     \end{itemize}
\end{itemize}
 We will show that if there is a collection that contains $l$ cubes with positive $\mu$ measure, then the same collection also contains $l+1$ cubes with positive $\mu$ measure. Due to Lemma \ref{PHprinc} we know that case $(2)$ is proved for $l=1,2,3$. Assume that case $(2)$ holds for $l$, such that $l \geq 3$.
 \begin{cla} \label{claim}
 There is a collection in which there are at least $l+1$ cubes with positive measure.
 \end{cla}
 Suppose that at every iteration we cannot find $l+1$ cubes with positive measure. If we fail to find an $l+1-$th cube at the $n-$th iteration, we obtain a cube, say $Q^{(n)}_{1k}$, of side-length $\frac{1}{6^{n}}$ for which $\mu(Q^{(n)}_{1k}) \geq \frac{1}{l^{n}3^{n(d-1)}2^{2n}}$. By the Frostman measure condition we have $ \frac{1}{l^{n}3^{n(d-1)}2^{2n}} \leq \mu(Q^{(n)}_{1k}) \leq C_{\mu} \frac{1}{6^{ns}}$, from which we obtain $ n \leq \frac{\log_2(C_{\mu})}{(1+\log_2(3))s-(d-1)\log_2(3) - \log_2(l)-2}$ for every $n$ which is a contradiction. \qedd

\bigskip

\end{document}